\documentclass[12pt]{amsart}
\usepackage{amscd, amsfonts, amssymb}
\usepackage{fullpage}
\usepackage{latexsym}
\usepackage{verbatim}
\usepackage{enumerate}
\usepackage[colorlinks=true,citecolor=red,linkcolor=blue]{hyperref}
\usepackage{tikz}
\usetikzlibrary{arrows,matrix,decorations,positioning, backgrounds,fit}

\newcommand\ra{\rightarrow}

\DeclareMathOperator{\Ad}{Ad}

\DeclareMathOperator{\Char}{char}

\DeclareMathOperator{\HH}{H}
\DeclareMathOperator{\Z}{Z}
\DeclareMathOperator{\GL}{GL}

\DeclareMathOperator{\HHom}{Hom}

\numberwithin{equation}{section}

\newtheorem{thm}[equation]{Theorem}

\theoremstyle{definition}

\theoremstyle{remark}

\theoremstyle{remark}
\newtheorem{rems}[equation]{Remarks}

\newtheorem{qn}[equation]{Question}

\newcommand{\ovl}{\overline}


\title[On a question of K\"ulshammer]
{On a question of K\"ulshammer for representations of finite groups in reductive groups}

\author[M.\  Bate]{Michael Bate}
\address
{Department of Mathematics,
University of York,
York YO10 5DD,
United Kingdom}
\email{michael.bate@york.ac.uk}

\author[B.\ Martin]{Benjamin Martin}
\address
{Department of Mathematics,
University of Aberdeen,
King's College,
Fraser Noble Building,
Aberdeen AB24 3UE,
United Kingdom}
\email{b.martin@abdn.ac.uk}

\author[G. R\"ohrle]{Gerhard R\"ohrle}
\address
{Fakult\"at f\"ur Mathematik,
Ruhr-Universit\"at Bochum,
D-44780 Bochum, Germany}
\email{gerhard.roehrle@rub.de}

\dedicatory{To Burkhard K\"ulshammer on his sixtieth birthday}

\begin{document}

\begin{abstract}
Let $G$ be a simple algebraic group of type $G_2$ over an algebraically closed field of characteristic $2$.  
We give an example of a finite group $\Gamma$ with Sylow $2$-subgroup $\Gamma_2$ and an infinite family of pairwise 
non-conjugate homomorphisms $\rho\colon \Gamma\ra G$ whose restrictions to $\Gamma_2$ are all conjugate.  
This answers a question of Burkhard K\"ulshammer from 1995.  We also give an action of $\Gamma$ on a connected unipotent group $V$ such that the map of 1-cohomologies $\HH^1(\Gamma,V)\ra \HH^1(\Gamma_p,V)$ induced by restriction of 1-cocycles has an infinite fibre.
\end{abstract}

\maketitle

\section{Introduction}
\label{sec:intro}

Let $k$ be an algebraically closed field and let $\Gamma$ be a finite group.  
By a {\em representation} of $\Gamma$ in a linear algebraic group $H$ over $k$, 
we mean a group homomorphism from $\Gamma$ to $H$.  We denote by 
$\HHom(\Gamma,H)$ the set of representations $\rho$ of $\Gamma$ in $H$; 
this has the natural structure of an affine variety over $k$ (see, e.g., \cite[II.2]{slodowy}).  
The group $H$ acts on $\HHom(\Gamma,H)$ 
by conjugation and we call the orbits $H\cdot \rho$ {\em conjugacy classes}.

If either $\Char(k)= 0$ or $\Char(k)= p> 0$ and $|\Gamma|$ is coprime 
to $p$, then every representation of $\Gamma$ in $\GL_n(k)$ is 
completely reducible and $\HHom(\Gamma,\GL_n(k))$ is a 
finite union of conjugacy classes, by Maschke's Theorem.  
Now suppose that $\Char(k)= p> 0$ and $p$ divides $|\Gamma|$.  
It is no longer 
true that $\HHom(\Gamma,\GL_n(k))$ is a finite union of 
conjugacy classes---for example, this fails even for $n= 2$ and 
$\Gamma= C_p\times C_p$ 
(cf.\ the last paragraph of the proof of Theorem~\ref{thm:main} below).  Let $\Gamma_p$ be a Sylow $p$-subgroup of $\Gamma$.  
It is natural to ask instead whether representations of $\Gamma$ 
are controlled by their restrictions to $\Gamma_p$.  
Burkhard K\"ulshammer raised the following question in 1995 in \cite[Sec.~2]{kuls} (see also \cite[I.5]{slodowy}).

\begin{qn}
\label{qn:K}
 Let $G$ be a linear algebraic group and let $\sigma\in \HHom(\Gamma_p,G)$.  
Are there only finitely many conjugacy classes of representations 
$\rho\in \HHom(\Gamma,G)$ such that $\rho|_{\Gamma_p}$ is conjugate to $\sigma$?
\end{qn}

Straightforward representation-theoretic arguments show that the 
answer is yes if $G= \GL_n(k)$ (see \cite[Sec.~2]{kuls}).  
On the other hand, an example of Cram with $p= 2$ shows that the answer is no in general if we allow $G$ to be non-connected and non-reductive \cite{cram}.

For the rest of this paper, we assume $G$ is connected and reductive. 
Slodowy proved that the answer to Question~\ref{qn:K} is yes under some extra hypotheses \cite{slodowy}; we briefly summarise his results.  If one embeds $G$ in some $\GL_n(k)$, 
then $\HHom(\Gamma,G)$ embeds in $\HHom(\Gamma,\GL_n(k))$.  
Given $\rho\in \HHom(\Gamma,G)$, the set $(\GL_n(k)\cdot \rho)\cap \HHom(\Gamma,G)$ 
splits into a union of $G$-conjugacy classes; in the first part of his paper, Slodowy applies a beautiful geometric argument due to Richardson \cite{rich2} 
to show that this union is finite when $p$ is good for $G$, which allows one to deduce a positive answer 
to Question~\ref{qn:K} for $G$ from the positive answer for $\GL_n(k)$ \cite[I.5, Thm.\ 3]{slodowy}.  

The second part of Slodowy's paper gives a different criterion for Question~\ref{qn:K} 
to have positive answer: he shows that if $\sigma(\Gamma_p)$ has reduced centralizer 
in $G$ then there are only finitely many conjugacy classes of representations $\rho\in \HHom(\Gamma,G)$ 
such that $\rho|_{\Gamma_p}$ is conjugate to $\sigma$ \cite[II.4, Cor.\ 1]{slodowy}.  An important ingredient in this proof, which dates back to work of Andr\'e Weil, 
is that one can interpret elements of the tangent space to $\HHom(\Gamma,G)$ at $\rho$
as elements of the space of 1-cocycles $\Z^1(\Gamma,{\mathfrak g})$, where ${\mathfrak g}$ 
denotes the Lie algebra of $G$ and $\Gamma$ acts on ${\mathfrak g}$ by $\gamma\cdot X=\Ad(\rho(\gamma))(X)$.
In fact, Slodowy proved a more general finiteness criterion 
in terms of the ``inseparability defects'' of $\rho$ and $\rho|_{\Gamma_p}$ \cite[II.4, Thm.\ 2]{slodowy}\footnote{This result actually holds for non-reductive $G$ as well.}.  The case of arbitrary connected reductive $G$ was, however, still left open.  

In this note we show that the answer to Question~\ref{qn:K} is no in general for connected reductive $G$.  
We prove the following result.

\begin{thm}
\label{thm:main}
 Suppose $G$ is a simple algebraic group of type $G_2$ and $\Char(k)= 2$.  Let $q> 3$ 
be odd, let $D_{2q}$ denote the dihedral group of order $2q$, let 
$\Gamma= D_{2q}\times C_2= \langle r,s,z\mid r^q= s^2= z^2= 1, 
srs^{-1}= r^{-1}, [r,z]= [s,z]= 1\rangle$ 
and let $\Gamma_2= \langle s,z\rangle$ (a Sylow $2$-subgroup of $\Gamma$).  
Then there exist representations $\rho_a\in \HHom(\Gamma,G)$ for all $a\in k$ 
such that the $\rho_a$ are pairwise non-conjugate and the restrictions $\rho_a|_{\Gamma_2}$ are conjugate for all $a\in k$.
\end{thm}

In \cite[Sec.~7]{BMRT} the authors and Tange
constructed families of finite subgroups of $G= G_2$ in characteristic $2$ with unusual properties (note, for example, that \cite[Ex.\ 7.15]{BMRT} shows that Richardson's argument can fail in bad characteristic).  
Our proof of Theorem~\ref{thm:main} involves a modification of this construction.

Our results can be interpreted in the language of nonabelian 1-cohomology (see Section~\ref{sec:nonabelian}).  Let $\Gamma$ act by group automorphisms on a unipotent group $V$.  One can form the 1-cohomology $\HH^1(\Gamma,V)$, and the inclusion of $\Gamma_p$ in $\Gamma$ gives a map $\Theta$ from $\HH^1(\Gamma,V)$ to $\HH^1(\Gamma_p,V)$ induced by restriction of 1-cocycles.

\begin{thm}
\label{thm:H1fibre}
 Let $p= 2$, let $q> 3$ be odd and let $\Gamma= D_{2q}\times C_2$.  There is an action of $\Gamma$ on a connected unipotent group $V$ such that the map $\Theta$ has an infinite fibre.
\end{thm}

\noindent This is in sharp contrast to the case when $V$ is abelian: standard results from abelian cohomology (cf.\ \cite[III, Prop.\ 10.4]{brown}) show that if $V$ is an abelian unipotent group (e.g., a finite-dimensional vector space over $k$) on which $\Gamma$ acts by group automorphisms then $\Theta$ is injective.  In fact, Slodowy uses precisely this result in the special case when $V$ is the $\Gamma$-module ${\mathfrak g}$ on the way to proving \cite[II.4, Thm.\ 2]{slodowy} (see \cite[II.4, Lem.]{slodowy}).

Lond gave a different example with $\Theta$ having an infinite fibre \cite[Ex.\ 4.1]{lond}, using the example of Cram discussed above.  In our case, the group $V$ is the unipotent radical of a parabolic subgroup $P$ of a simple group $G$ of type $G_2$, and $\Gamma$ acts on $V$ by conjugation, via a homomorphism $\sigma\colon \Gamma\ra P$.  Theorem~\ref{thm:H1fibre} follows quickly from the construction in Section~\ref{sec:proof} (see Section~\ref{sec:nonabelian}).

\section{Proof of Theorem~\ref{thm:main}}
\label{sec:proof}

Until the end of this section we take $G$ to be a simple algebraic group of type $G_2$ and $\Char(k)$ to be $2$.  
We recall some notation from \cite[Sec.\ 7]{BMRT}. 
The positive roots of $G$ with respect to a fixed maximal torus 
$T$ and a fixed Borel subgroup containing $T$ are 
$\alpha$ (short), $\beta$ (long), $\alpha+ \beta$, $2\alpha+ \beta$, $3\alpha+ \beta$ and 
$\omega := 3\alpha+ 2\beta$.  
Given a root $\delta$, we denote the corresponding root group by $U_\delta$ 
and coroot by $\delta^\vee$.  
We fix a group isomorphism $\kappa_\delta \colon k\ra U_\delta$.  
We write $G_\delta$ for $\langle U_\delta\cup U_{-\delta}\rangle$ 
and we set $s_\delta = \kappa_\delta(1) \kappa_{-\delta}(1) \kappa_\delta(1)$; then $s_\delta$ 
represents the reflection corresponding to $\delta$ in the Weyl group of $G$ 
(since $\Char(k)= 2$, $s_\delta$ has order $2$).

Fix $t\in \alpha^\vee(k^*)$ such that $|t|= q$.  For $a\in k$, define $\rho_a\in \HHom(\Gamma,G)$ by
\[
\rho_a(r)= t, \quad \rho_a(s)= s_\alpha \kappa_{\omega}(a), \quad \rho_a(z)= \kappa_{\omega}(1). 
\]
It is easily checked that this is well-defined (note that $[G_\alpha, G_{\omega}]= 1$).  
Set $u(x)= \kappa_\beta(x) \kappa_{3\alpha+ \beta}(x)$ for $x\in k$.  
Then $u(x)$ commutes with $U_{\omega}$ and 
$u(x)s_\alpha u(x)^{-1}= s_\alpha \kappa_{\omega}(x^2)$ 
(see the first paragraph of \cite[p.~4307]{BMRT}).  
It follows that $u(\sqrt{a})\cdot \left(\rho_0|_{\Gamma_2}\right)= \rho_a|_{\Gamma_2}$.

To complete the proof of Theorem~\ref{thm:main}, we now need to show that the 
$\rho_a$ are pairwise non-conjugate.  Let $a,b\in k$
and suppose $g\cdot \rho_a= \rho_b$ for some $g\in G$.  
Then $g\in C_G(t)$.  It follows from \cite[(7.1) and (7.2)]{BMRT} that $C_G(t)= TG_{\omega}$ (this is where we need our assumption that $q> 3$; cf.\ \cite[(7.7)]{BMRT}).  
So write $g= hm$ with $h\in T$ and $m\in G_{\omega}$.  
We have $(hm)s_\alpha \kappa_{\omega}(a)(hm)^{-1}= s_\alpha \kappa_{\omega}(b)$, so $hs_\alpha h^{-1} (hm)\kappa_{\omega}(a)(hm)^{-1}= s_\alpha \kappa_{\omega}(b)$ since $m$ commutes with $s_\alpha$.  Now $G_\alpha\cap G_\omega= 1$ (see the paragraph following \cite[(7.8)]{BMRT}), so the condition $hs_\alpha h^{-1} (hm)\kappa_{\omega}(a)(hm)^{-1}= s_\alpha \kappa_{\omega}(b)$ forces $h$ to commute with $s_\alpha$, as $hs_\alpha h^{-1}\in G_\alpha$ and $(hm)\kappa_{\omega}(a)(hm)^{-1}\in G_\omega$.  A simple calculation now shows that $h\in {\rm ker}(\alpha)\subseteq G_{\omega}$.  
Hence $g\in G_{\omega}$.  
But $G_{\omega}$ is a simple group of type $A_1$, so the pair 
$(\kappa_{\omega}(a), \kappa_{\omega}(1))$ is not 
$G_{\omega}$-conjugate to the pair $(\kappa_{\omega}(b), \kappa_{\omega}(1))$ unless $a= b$.  
We conclude that $\rho_a$ and $\rho_b$ are not conjugate if $a\neq b$, as required.

\begin{rems}
 (i). Choose an embedding $i$ of $G$ in some $\GL_n(k)$.  
Then the representations $i\circ \rho_a$ of $\Gamma$ in $\GL_n(k)$ 
fall into finitely many $\GL_n(k)$-conjugacy classes, 
since Question~\ref{qn:K} has positive answer for $\GL_n(k)$.  
Hence there exists $a\in k$ such that $\left(\GL_n(k)\cdot \rho_a\right)\cap \HHom(\Gamma,G)$ 
is an infinite union of $G$-conjugacy classes.  
This gives another example of the phenomenon in \cite[Ex.\ 7.15]{BMRT} discussed above.
 
 (ii). It follows from Slodowy's result \cite[II.4, Thm.\ 2]{slodowy} discussed above that $\rho_a$ has greater inseparability defect than $\rho_a|_{\Gamma_2}$ for at least one $a\in k$.  In fact, it can be shown using the calculations in \cite[Sec.\ 7]{BMRT} that if $a\neq 0$ then $\rho_a$ has inseparability defect 1 and $\rho_a|_{\Gamma_2}$ has inseparability defect 5.  This answers a question of Slodowy \cite[II.4, Rem.\ 2]{slodowy}.
\end{rems}

We do not know of any analogous examples in odd characteristic; recall from the discussion in Section~\ref{sec:intro} that if such an example exists then $p$ must be bad for $G$.  Our construction is closely related to the construction of a certain triple $(G,M,H)$ in \cite[Sec.~7]{BMRT}, where $G= G_2$, $M$ is a reductive subgroup of $G$ and $H$ is a finite subgroup of $M$.  We guess that further examples can be obtained from other triples $(G,M,H)$ with similar properties, but we leave this for future work.  The mechanism for producing these triples works only in characteristic 2 (see the paragraph following \cite[Rem.\ 1.6]{uchiyama2}).  Uchiyama found triples $(G,M,H)$ for $G$ of type $E_7$ \cite[Sec.\ 3]{uchiyama2}, and showed that the construction fails for several cases involving groups of rank at most 6, including $A_3$, $A_4$, $B_3$ and $E_6$ \cite[Thm.\ 3.1.1, Ch.\ 4]{uchiyama1}.

It seems an interesting problem to find examples like that of Cram \cite{cram} but in odd characteristic, where we allow $G$ to be non-reductive.

\section{Nonabelian 1-cohomology}
\label{sec:nonabelian}

Another approach to K\"ulshammer's problem is via the $1$-cohomology 
of the unipotent radical $R_u(P)$, where $P$ is a proper parabolic subgroup of $G$.  
Here is a brief explanation.  Recall that a closed subgroup $M$ of $G$ is said to be 
{\em $G$-completely reducible} if whenever $M$ is contained 
in a parabolic subgroup $P$ of $G$, $M$ is contained in some Levi subgroup of $P$ 
\cite{serre2}, \cite{serre1}.  As a special case, we say that $M$ is {\em $G$-irreducible} 
if $M$ is not contained in any proper parabolic subgroup of $G$ at all.  
We say that $\rho\in \HHom(\Gamma,G)$ is $G$-completely reducible 
(resp., $G$-irreducible) if its image is.

Although in general $\HHom(\Gamma,G)$ is an infinite union of conjugacy classes 
for reductive $G$, it was proved in \cite[Cor.\ 3.8]{BMR} that there are only 
finitely many conjugacy classes of representations that are $G$-completely reducible.  
This generalizes the classical result that a finite group admits only finitely many 
completely reducible $n$-dimensional representations in any characteristic.  
Moreover, it follows from \cite[Cor.\ 3.7]{BMR} that the conjugacy classes of 
$G$-completely reducible representations of $\Gamma$ in $G$ are precisely 
the conjugacy classes that are Zariski-closed subsets of $\HHom(\Gamma,G)$.  
Given $\rho\in \HHom(\Gamma,G)$, choose a minimal parabolic subgroup 
$P$ of $G$ with $\rho(\Gamma)\subseteq P$.  Let $L$ be a Levi subgroup of 
$P$ and let $\pi\colon P\ra L$ be the canonical projection.  
It follows from \cite[Cor.\ 3.5]{BMR} that $\sigma:= \pi\circ\rho\in \HHom(\Gamma,L)$ 
is $L$-irreducible and $G$-completely reducible.  
Conversely, given $G$-irreducible $\sigma\in \HHom(\Gamma,G)$, 
we can consider the set $C_\sigma$ of all $\rho\in \HHom(\Gamma,P)$ 
such that $\pi\circ \rho= \sigma$.  
By the result described in the first sentence of this paragraph, there are only finitely many possibilities 
for $(P, L, \sigma)$ up to $G$-conjugacy.  
Hence if $C\subseteq \HHom(\Gamma,G)$ is an infinite union of 
$G$-conjugacy classes then for some triple $(P, L, \sigma)$, 
$C_\sigma$ must meet infinitely many $G$-conjugacy classes in $C$.  
Thus we have reduced the ``global'' problem of considering all representations 
into $G$ to the ``local'' problem of considering all representations into a fixed proper parabolic subgroup $P$.

Next we study the structure of $C_\sigma$ for fixed $(P,L,\sigma)$.  
Let $V= R_u(P)$.  Given $\rho\in C_\sigma$, there is a unique function 
$\theta_\rho\colon \Gamma\ra V$ 
defined by $\rho(\gamma)= \theta_\rho(\gamma)\sigma(\gamma)$.  
It is easily checked that $\theta_\rho$ satisfies the $1$-cocycle relation 
$\theta_\rho(\gamma_1\gamma_2)= \theta_\rho(\gamma_1)(\gamma_1\cdot \theta_\rho(\gamma_2))$, 
where $\Gamma$ acts on $V$ by $\gamma\cdot v= \sigma(\gamma)v\sigma(\gamma)^{-1}$.  
The converse is also true, so we have a bijection between $C_\sigma$ 
and the space of $1$-cocycles $\Z^1(\Gamma,\sigma,V)$.  
A simple calculation shows that $\rho, \mu\in C_\sigma$ are $V$-conjugate 
if and only if the images $\ovl{\theta_\rho}$ of $\theta_\rho$ and 
$\ovl{\theta_\mu}$ of $\theta_\mu$ in $\HH^1(\Gamma,\sigma,V)$ are equal.  
Thus we have an interpretation of $V$-conjugacy classes in $C_\sigma$ in terms of $1$-cohomology (cf.\ the proof of \cite[I.5, Lem.\ 1]{slodowy}).

This idea has been used in a slightly different context to study embeddings of reductive 
algebraic groups inside simple algebraic groups \cite{liebeckseitz0}, \cite{stewartG2}, \cite{stewartF4}, \cite{lond}.  
In our case we have an extra ingredient arising from restriction of representations.  
The restriction map from $\HHom(\Gamma,G)$ to $\HHom(\Gamma_p,G)$ 
maps $C_\sigma$ to $C_{\sigma|_{\Gamma_p}}$.  
Restriction of cocycles gives a map from $\Z^1(\Gamma,\sigma,V)$ to $\Z^1(\Gamma_p,\sigma|_{\Gamma_p},V)$ 
which is compatible with the correspondence between representations and $1$-cocycles, 
and this descends to give a map $\Theta$ from $\HH^1(\Gamma,\sigma,V)$ to $\HH^1(\Gamma_p,\sigma|_{\Gamma_p},V)$.  
See \cite[Ch.\ 3--4]{lond} for a fuller explanation.

Now we recast our example in this language.  
Let $G$, $k$, $\Gamma$, $\Gamma_2$ and the $\rho_a$  be as in Section~\ref{sec:proof}.  
Set $P= P_\alpha$, $L= L_\alpha$ and $V= R_u(P_\alpha)$, 
and define $\sigma\in \HHom(\Gamma,L)$ by $\sigma(r)= t$, $\sigma(s)= s_\alpha$ 
and $\sigma(z)= 1$.  Then $\sigma$ is $L$-irreducible and every $\rho_a$ belongs to $C_\sigma$.  
Let $\theta_a\in \Z^1(\Gamma,\sigma,V)$ and $\theta'_a\in \Z^1(\Gamma_2,\sigma|_{\Gamma_2},V)$ 
be the $1$-cocycles corresponding to $\rho_a$ and $\rho_a|_{\Gamma_2}$, respectively.  
The calculations in Section~\ref{sec:proof} show that the $\rho_a|_{\Gamma_2}$ 
are pairwise $V$-conjugate, so the $1$-cohomology classes $\ovl{\theta'_a}\in \HH^1(\Gamma_2,\sigma|_{\Gamma_2},V)$ 
are equal for all $a\in k$.  In contrast, no two of the $\rho_a$ are $V$-conjugate 
(since no two are $G$-conjugate), so the $1$-cohomology classes $\ovl{\theta_a}\in \HH^1(\Gamma,\sigma,V)$ are all different.  
Thus we have an example where the map $\Theta$ from $\HH^1(\Gamma,\sigma,V)$ to $\HH^1(\Gamma_2,\sigma|_{\Gamma_2},V)$ 
has an infinite fibre (cf.\  \cite[Ex.\ 4.1]{lond}).

We do not know of any analogous examples in odd characteristic; cf.\ the discussion at the end of Section~\ref{sec:proof}.


\bigskip
{\bf Acknowledgments}:
The authors acknowledge the financial support of EPSRC Grant EP/L005328/1, 
 Marsden Grants UOC1009 and UOA1021, and
the DFG-priority programme SPP1388 ``Representation Theory''.  We are grateful to the referee for helpful suggestions.


\end{document}